\documentclass[12pt]{article}
\usepackage{amssymb, amsfonts, amsmath, amsthm}

\def\3{\subset }
\def\4{\subseteq }

\def\0{\leqno}

\def\barr{\begin{array}}
\def\earr{\end{array}}
\def\dd{\displaystyle}

\def\Z{{\rlap{$\kern2pt{\rm Z}$}{\rm Z}\,}}

\title{A characterization of ZM-groups}
\author{Marius T\u arn\u auceanu}
\date{April 11, 2017}

\begin{document}

\maketitle

\begin{abstract}
In this short note we give a characterization of ZM-groups that uses the functions defined and studied in \cite{3,4}. This leads to a proof of Conjecture 6 in \cite{4}.
\end{abstract}



\section{Introduction}

Given a finite group $G$, we consider the functions
$$\varphi(G)=|\{a\in G \mid o(a)=\exp(G)\}| \mbox{ and } S(G)=\dd\sum_{H\leq\, G}\varphi(H).$$The class of finite groups $G$ for which $S(G)=|G|$ has been partially determined in \cite{4}. We are now able to complete this study by proving the following result.

\bigskip\noindent{\bf Theorem 1.1.} {\it Let $G$ be a finite group. Then $S(G)\geq |G|$, and we have equality if and only if $G$ is a ZM-group,
i.e. a group with all Sylow subgroups cyclic.}
\bigskip

In particular, since a ZM-group is nilpotent if and only if it is cyclic, we obtain the result in Conjecture 6 of \cite{4}.

\bigskip\noindent{\bf Corollary 1.2.} {\it Let $G$ be a finite nilpotent group. Then $S(G)\geq |G|$, and we have equality if and only if $G$ is cyclic.}

\section{Proof of Theorem 1.1}

Let $C(G)$ be the poset of cyclic subgroups of $G$. For every divisor $d$ of $|G|$ we denote by $n_d$ the number of cyclic subgroups of order $d$ of $G$ and by $n'_d$ the number of elements of order $d$ in $G$. Then we have $$n'_d=n_d\varphi(d)$$because a cyclic subgroup of order $d$ contains $\varphi(d)$ elements of order $d$. One obtains $$S(G)\geq\hspace{-1mm}\dd\sum_{H\in C(G)}\hspace{-1mm}\varphi(H)=\hspace{-1mm}\dd\sum_{H\in C(G)}\hspace{-1mm}\varphi(|H|)=\dd\sum_{d\,\mid\,|G|}n_d\varphi(d)=\dd\sum_{d\,\mid\,|G|}n'_d=|G|,$$as desired.

Assume now that $S(G)=|G|$. Then $\varphi(H)=0$ for all non-cyclic subgroups $H$ of $G$. Since for any $p$-group $P$ we have $\varphi(P)\neq 0$ (see \cite{3}), it follows that all Sylow subgroups of $G$ are cyclic, that is $G$ is a ZM-group.

Conversely, assume that $G$ is a ZM-group. By \cite{2} such a group is of type
$${\rm ZM}(m,n,r)=\langle a, b \mid a^m = b^n = 1, \hspace{1mm}b^{-1} a b = a^r\rangle,$$where the triple $(m,n,r)$ satisfies the conditions
$${\rm gcd}(m,n)={\rm gcd}(m,r-1)=1 \mbox{ and } r^n\equiv 1 \hspace{1mm}({\rm mod}\hspace{1mm}m).$$It is clear that $|{\rm ZM}(m,n,r)|=mn$.
The subgroups of ${\rm ZM}(m,n,r)$ have been completely described in \cite{1}. Set
$$L=\left\{(m_1,n_1,s)\in\mathbb{N}^3 \hspace{1mm}\mid\hspace{1mm} m_1|m,\hspace{1mm} n_1|n,\hspace{1mm} s<m_1,\hspace{1mm} m_1|s\frac{r^n-1}{r^{n_1}-1}\right\}.$$Then there is a bijection between $L$ and the subgroup lattice $L({\rm ZM}(m,n,r))$ of ${\rm ZM}(m,n,r)$, namely the function that maps a triple $(m_1,n_1,s)\in L$ into the subgroup $H_{(m_1,n_1,s)}$ defined by
$$H_{(m_1,n_1,s)}=\bigcup_{k=1}^{\frac{n}{n_1}}\alpha(n_1,s)^k\langle a^{m_1}\rangle=\langle a^{m_1},\alpha(n_1, s)\rangle,$$where $\alpha(x, y)=b^xa^y$, for all $0\leq x<n$ and $0\leq y<m$. We can easily check that
$$|H_{(m_1,n_1,s)}|=\exp(H_{(m_1,n_1,s)})=\frac{mn}{m_1n_1}\,,$$and so $$\varphi(H_{(m_1,n_1,s)})\neq 0 \mbox{ if and only if } H_{(m_1,n_1,s)} \mbox{ is cyclic}.$$This shows that
$$S(G)=\hspace{-1mm}\dd\sum_{H\in C(G)}\hspace{-1mm}\varphi(H)=|G|,$$completing the proof.
\hfill\rule{1,5mm}{1,5mm}

\vspace*{5ex}\small

\hfill
\begin{minipage}[t]{5cm}
Marius T\u arn\u auceanu \\
Faculty of  Mathematics \\
``Al.I. Cuza'' University \\
Ia\c si, Romania \\
e-mail: {\tt tarnauc@uaic.ro}
\end{minipage}

\end{document}